\documentclass[a4paper,10pt]{article}
\usepackage{amsmath,amsthm,amsfonts}
\usepackage{graphicx}

\theoremstyle{plain}
\newtheorem{thm}{Theorem}

\theoremstyle{remark}

\begin{document}

\title{Reductions of $(v_3)$ configurations}

\author{Marko Boben\thanks{Address: Marko Boben, Univerza v Ljubljani,
Fakulteta za ra\v{c}unalni\v{s}tvo in informatiko,
Tr\v{z}a\v{s}ka 25, 1000 Ljubljana, Slovenia;
E-mail address: \texttt{Marko.Boben@fri.uni-lj.si}}\\
University of Ljubljana, FRI, IMFM\\
Slovenia}

\date{\today}

\maketitle

\begin{abstract}
Cubic bipartite graphs with girth at least $6$ correspond to
symmetric combinatorial $(v_3)$ configurations. In 1887 V.\ Martinetti
described a simple reduction method which enables one to reduce
each combinatorial $(v_3)$ configuration to one from the infinite set
of so-called \emph{irreducible} configurations.
The aim of this paper is to show that a slightly extended set of reductions
enables one to reduce each combinatorial $(v_3)$ configuration either to the Fano
configuration or to the Pappus configuration.
\end{abstract}

\begin{flushleft}
\textbf{Key words:} cubic graphs, incidence structures, configurations.\\
\textbf{Math.\ Subj.\ Class.\ (2000):} 05C62, 05C85, 68R05, 68R10.
\end{flushleft}

\section{Introduction}

Although the results and methods presented in this article come out from
the graph theory, it is maybe worth mentioning that the problem itself
originates from incidence structures named combinatorial configurations
and dates back to 1887~\cite{Martinetti}.

A \emph{symmetric combinatorial configuration} of type $(v_r)$ is an incidence
structure of $v$ points and $v$ blocks (called lines) such that each line
contains $r$ points, each point belongs to $r$ lines and each two different
points are contained in at most one line.

One way to relate graphs and configurations is via their incidence graphs
(or Levi graphs~\cite{PiRa}). An \emph{incidence graph} or \emph{Levi graph} of the
$(v_r)$ configuration $\mathcal{C}$ is a bipartite $r$-regular graph on $v$
\emph{black} vertices representing the points of $\mathcal{C}$, $v$ \emph{white}
vertices representing the lines of $\mathcal{C}$, and with an edge joining two
vertices if and only if the corresponding point and line are incident in $\mathcal{C}$.

The property that two different points of a configuration are contained in at most
one line implies that Levi graphs have girth (the length of the shortest cycle)
at least $6$. Conversely, each bipartite $r$-regular graph with girth at least $6$ and
with chosen black-white coloring determines precisely one $(v_r)$ configuration.
Note that if the coloring is not chosen, such graphs determine a pair of
\emph{dual configurations}, i.e.\ configurations with the role of their points
and lines interchanged. More about graphs and configurations can be found
in~\cite{BBP,Coxeter,GroppCG}.

In 1887 V.\ Martinetti described the following construction method for symmetric $(v_3)$
configurations~\cite{GruNotes,Martinetti}. Suppose that in the given $(v_3)$ configuration there
exist two parallel (non-intersecting) lines $\{a_1,a_2,a_3\}$ and $\{b_1,b_2,b_3\}$ such
that the points $a_1$ and $b_1$ are not on a common line. By removing these two lines,
adding one new point $c$ and three new lines $\{c,a_2,a_3\}$,
$\{c,b_2,b_3\}$, $\{c,a_1,b_1\}$ we obtain a $((v+1)_3)$ configuration.
It is not possible to obtain every $(v_3)$ configuration from some $((v-1)_3)$
configuration by using this method. We call $(v_3)$ configurations which can not be
constructed in this way from a smaller one \emph{A-irreducible configurations}
and the others \emph{A-reducible configurations}. We use the prefix A before the
term (ir)reducible to distinguish this type of reduction from a more general
one defined later in the paper.

Let us consider the Martinetti method on Levi graphs of $(v_3)$ configurations.
For the sake of simplicity we will use, as in~\cite{MB}, the notion \emph{$(v_3)$ graph}
instead of Levi graph of $(v_3)$ configuration, i.e.\ $(v_3)$ graph is a bipartite cubic graph with
girth $\geq 6$. We define A-reducible and A-irreducible $(v_3)$ graphs corresponding
to A-reducible and A-irreducible configurations respectively as follows.
A $(v_3)$ graph $G$ is \emph{A-reducible} if there exists an edge $uv \in E(G)$ such that
$(G - \{u,v\}) + x_1 y_1 + x_2 y_2$ or $(G - \{u,v\}) + x_1 y_2 + x_2 y_1$
is again a $(v_3)$ graph, where $x_1,x_2$ are the neighbors of $u$ (different from $v$)
and $y_1,y_2$ are the neighbors of $v$ (different from $u$).
Otherwise we call a $(v_3)$ graph \emph{A-irreducible}.

It turns out that there are four infinite families of connected A-irreducible graphs,
which we present in the next section, plus the Pappus graph,
i.e.\ the Levi graph of the Pappus configuration, see Figure~\ref{fig:pappusg}.
Note that Martinetti in his paper did not describe all families of
irreducible graphs. This was not done until~\cite{MB}.

In this paper we show that a slightly more general definition of a $(v_3)$ graph reduction
(and hence $(v_3)$ configuration reductions) leaves only two irreducible graphs:
the Heawood graph, the Levi graph of the Fano configuration (projective plane of order $2$),
and the already mentioned Pappus graph.

An effort in this direction was made in~2000 by H.\ G.\ Carstens, T.\ Dinski,
and E.\ Steffen~\cite{CDS}, who gave a set of quite complicated reductions to show that
they lead to the Heawood graph. But their claim is not entirely correct since
it turns out that their reductions leave out the Desargues graph too (the Levi graph of the
Desargues configuration)~\cite{erratum}.

\section{The main result}

Given a bipartite cubic graph $G$, black vertex $u \in V(G)$, and white vertex $v \in V(G)$,
let us say that we \emph{reduce $G$ by $u$ and $v$} if we remove $u$ and $v$ and
connect their neighbours in such way that $G$ remains cubic and bipartite.
If $u$ and $v$ are not connected then there are six ways to this, otherwise
there are only two.

If it is possible to reduce a $(v_3)$ graph $G$
by some pair of its vertices such that the reduced graph is again a $(v_3)$ graph,
then we say that $G$ is \emph{B-reducible}; otherwise we call it \emph{B-irreducible}.
We notice that A-reduction is only a special case of the above definition
where $u$ and $v$ are neighbouring vertices.

Our main result says that there are not many B-irreducible $(v_3)$ graphs.
\begin{thm} \label{thm:main}
There are only two connected B-irreducible bipartite cubic graphs with girth at least $6$,
the Heawood graph and the Pappus graph.
\end{thm}

To prove this theorem we show that our reduction applies to all A-irreducible
$(v_3)$ graphs but the Heawood and the Pappus graph.

It is possible that the reduction of a connected graph makes it disconnected.
So, the process of reducing a $(v_3)$ graph may finally result in a graph
whose connected components are only the Heawood and the Pappus graphs. But we may
state a stronger result.
\begin{thm} \label{thm:reduce}
Every connected $(v_3)$ graph can be reduced either to the Heawood graph or to the Pappus graph
by the sequence of connected $(v_3)$ graphs.
\end{thm}
\begin{proof}
Note that bipartite cubic graphs are 2-connected; this follows, for example,
from~\cite[Corollary 6.4.6]{Diestel}. Hence the removal of two vertices
produces at most three connected components in the general case and at most two
in the case of A-reduction. If the addition of the necessary edges leaves the graph
disconnected and this does not produce a $4$-cycle, then there also exists a choice
for the addition of the edges which gives a connected graph and does not produce a
$4$-cycle. (Such a choice exists since we add $3$ edges in the general case
and $2$ in the case of A-reduction.)
\end{proof}

\subsection{A-irreducible graphs}

Let us first list all A-irreducible graphs as they are presented in~\cite{MB}.

\paragraph{The first family of A-irreducible graphs.}

\begin{figure}[htb]
\centering
  \includegraphics{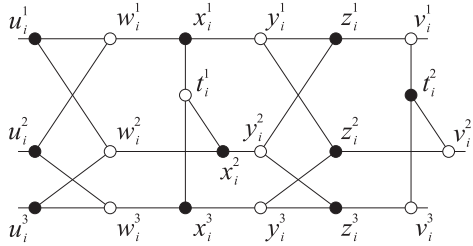}
  \caption{The segment graph $G_T$ defining the families $T_1(n)$, $T_2(n)$ and $T_3(n)$.}
  \label{fig:Tsegment}
\end{figure}

\begin{figure}[htb]
\centering
\begin{minipage}{0.3\textwidth}
  \centering \includegraphics{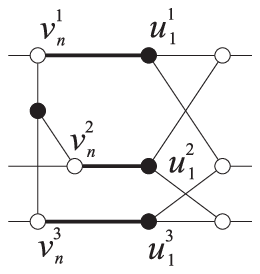}\\[2mm] (a)
\end{minipage}%
\begin{minipage}{0.3\textwidth}
  \centering \includegraphics{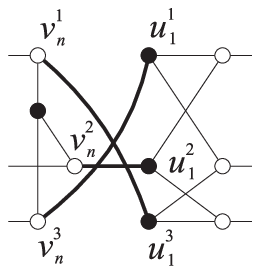}\\[2mm] (b)
\end{minipage}%
\begin{minipage}{0.3\textwidth}
  \centering \includegraphics{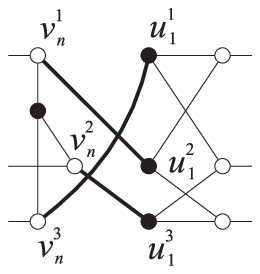}\\[2mm] (c)
\end{minipage}
\caption{The construction of the graphs $T_1(n)$ (a), $T_2(n)$ (b), and $T_3(n)$ (c)
from $T(n)$ by adding three edges (shown thick) joining the last and the first segment.}
\label{fig:TGraphs}
\end{figure}

The first type of A-irreducible graphs can be described as follows.
Let $T(n)$, $n \geq 1$, denote a graph on $20n$ vertices which is
an union of $n$ segments $G_T$ shown in Figure~\ref{fig:Tsegment} where
the $i$-th segment ($i \geq 2$) and the $(i-1)$-st segment are joined together
by the edges $v_{i-1}^1 u_i^1$, $v_{i-1}^2 u_i^2$, $v_{i-1}^3 u_i^3$.
We will use $T(n)$ in the following definitions.
Let $T_1(n)$ be the graph which is obtained from $T(n)$ by adding the edges
$u_1^1 v_n^1$, $u_1^2 v_n^2$, and $u_1^3 v_n^3$.
Let $T_2(n)$ be the graph obtained from $T(n)$ by adding the edges
$u^3_1 v^1_n$, $u^2_1 v^2_n$, and $u^1_1 v^3_n$.
And finally, let $T_3(n)$ be the graph obtained from $T(n)$ by adding the edges
$u^1_1 v^3_n$, $u^2_1 v^1_n$, and $u^3_1 v^2_n$. See Figure~\ref{fig:TGraphs}.

It turns out that for each fixed $n \geq 1$ the graphs $T_1(n)$, $T_2(n)$, $T_3(n)$
are pairwise non-isomorphic A-irreducible $(v_3)$ graphs.

\paragraph{The second family of A-irreducible graphs.}

\begin{figure} 
  \centering
  \begin{minipage}{0.5\textwidth}
    \centering
    \includegraphics{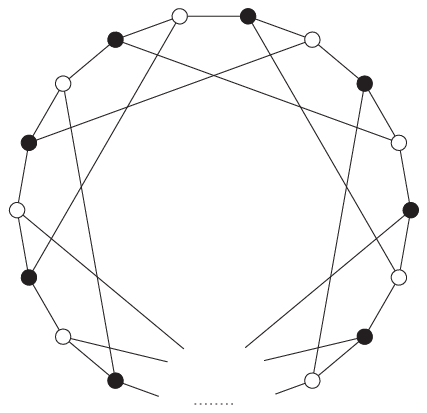}
  \end{minipage}%
  \begin{minipage}{0.5\textwidth}
    \centering
    \includegraphics{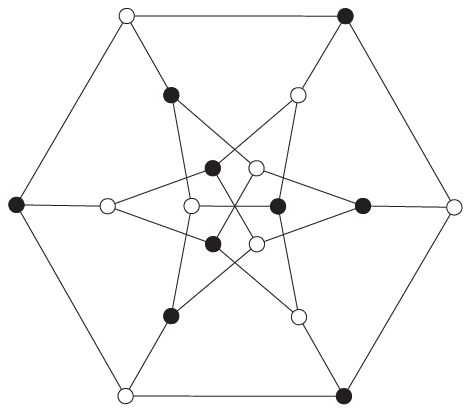}
  \end{minipage}
  \begin{minipage}[t]{0.5\textwidth}
    \centering
    \caption{Members of the second family of A-irreducible $(v_3)$ graphs are graphs $D(n)$
       with LCF notation $[5,-5]^n$, $n \geq 7$, \cite{Frucht}.}
    \label{fig:cgraph}
  \end{minipage}%
  \begin{minipage}[t]{0.5\textwidth}
    \centering
    \caption{The Pappus graph.}
    \label{fig:pappusg}
  \end{minipage}
\end{figure}

The second family of A-irreducible $(v_3)$ graphs consists of graphs $D(n)$ described by
LCF notation $[5, -5]^n$, $n \geq 7$, \cite{Frucht}. The notation $[5, -5]^n$ for $G$ means
that $G$ is constructed from the cycle of length $2n$ to which chords joining the vertices
at distance $5$ and $-5$ (in reverse direction) are attached, see Figure~\ref{fig:cgraph}.
The smallest graph of this family, $D(7)$, is the Heawood graph.
The members of this family are Levi graphs of the so-called cyclic $(v_3)$ configurations with base
line $\{0,1,3\}$.

\paragraph{The Pappus graph}

There exists another A-irreducible $(v_3)$ graph which is not a member of the families
described above. This is the Pappus graph, Levi graph of the Pappus configuration, see
Figure~\ref{fig:pappusg}.

The following theorem from~\cite{MB} states that the graphs mentioned above are
all A-irreducible graphs.
\begin{thm} \label{thm:Airreducible}
All connected A-irreducible $(v_3)$ graphs are
\begin{enumerate}
\item graphs $D(n)$, i.e.\ graphs with LCF notation $[5,-5]^n$, $n \geq 7$,
\item the graphs $T_1(n)$, $T_2(n)$, $T_3(n)$, $n \geq 1$, and
\item the Pappus graph.
\end{enumerate}
\end{thm}

\subsection{Proof of the main Theorem}

We conclude this section by the proof of Theorem~\ref{thm:main}.

\begin{proof}[Proof of Theorem~\ref{thm:main}]
Because of Theorem~\ref{thm:Airreducible} and the fact that A-reduction is a
special case of B-reduction, we have to take into consideration only A-irreducible
graphs.
The reduction of any $T(n)$ graph goes as follows: remove $w_1^2$
and $z_1^2$ and add edges $u_1^1 y_1^3$, $u_1^3 y_1^1$, and $v_1^2 x_1^2$,
see Figure~\ref{fig:reduction1}.

The reduction of any graph $D(n)$, $n \geq 8$, goes as follows: remove
$u_2^1$ and $u_2^4$ and add edges $v_1^1 u_2^3$, $u_2^2 v_2^2$, and $v_1^2 v_2^1$,
see Figure~\ref{fig:reduction2}.

\begin{figure} 
  \begin{minipage}{0.55\textwidth}
    \centering
    \includegraphics{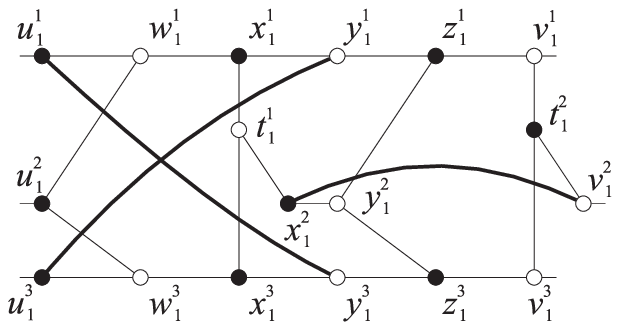}
  \end{minipage}%
  \begin{minipage}{0.45\textwidth}
    \centering
    \includegraphics{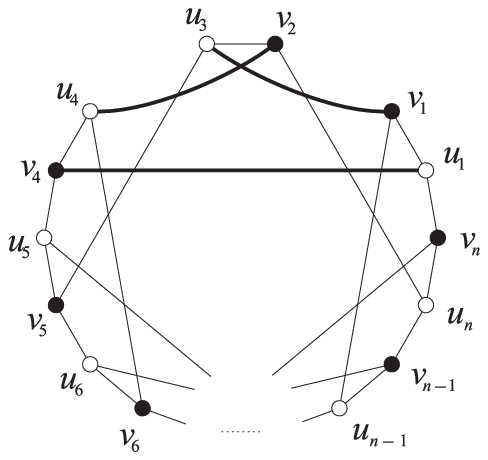}
  \end{minipage}
  \begin{minipage}[t]{0.55\textwidth}
    \centering
    \caption{Reduction of the first\newline A-irreducible family.}
    \label{fig:reduction1}
  \end{minipage} \hfil
  \begin{minipage}[t]{0.45\textwidth}
    \centering
    \caption{Reduction of the second\newline A-irreducible family.}
    \label{fig:reduction2}
  \end{minipage}
\end{figure}

In all cases it is easy to check that the girth of the reduced graph is still $6$.
In the case of the Heawood graph, $D(7)$, and the Pappus graph, any reduction
leads to a graph with girth $4$.
\end{proof}

\section{Conclusion}

In the discussion above we have proved, in the language of configurations,
that every combinatorial $(v_3)$ configuration can be reduced either to the Fano
configuration or to the Pappus configuration by a generalized Martinetti procedure.
In this procedure we always remove one line and one point and correct the structure
locally on the affected objects such that we again obtain a combinatorial configuration.

When speaking about combinatorial configurations we cannot avoid mentioning
\emph{geometric configurations}. A geometric $(v_r)$ configuration is the structure
of $v$ points and $v$ lines \emph{in the plane} such that each point belongs to $r$
lines and each line passes through $r$ points.
In the time when Martinetti wrote his article, the distinction between combinatorial
and geometric configurations was not strong, if it existed at all. One indication
that this was probably not the case is that he did not list the two smallest irreducible
configurations (the Fano $(7_3)$ configuration and the M\"{o}bius-Kantor $(8_3)$
configuration) since neither of them is a geometric configuration.

Hence we may ask the following question. For which combinatorial configurations that
can be realized as geometric configurations one can realize each of (at least one of)
the reduced configurations?
This is certainly not the case for the $(9_3)$ configurations, but do there exist
other ones?

Another question one may ask at the end is the following:
which cubic bipartite graphs of girth at least $g$ are irreducible if we do not allow
the girth of the reduced graph to fall under $g$? In Theorem~\ref{thm:main} we have found
graphs for $g = 6$, so the next step would be to find graphs for $g = 8$.
Obviously, one graph is the smallest trivalent graph of girth $8$ or \emph{the Tutte $8$-cage}
(which is also bipartite) on 30 vertices. The second one is the unique $(v_3)$ graph with girth
$8$ on $34$ vertices (since there is no such graph on $32$ vertices according to~\cite{BBP}).
There are also larger graphs as the generalized Petersen graph $G(18,5)$,
while $G(26,5)$ is reducible to a graph with girth $8$.

More about bipartite graphs of large girth and their connection to configurations can
be found in~\cite{PBMOA}.

\section*{Acknowledgements}

The research was supported in part by a grant from Ministrstvo za \v{s}olstvo,
znanost in \v{s}port Republike Slovenije, grant P1-0294. The author would also like
to thank to Tomo Pisanski for his encouragement and many fruitful suggestions during the
preparation of this paper.

\end{document}